\documentclass[12pt]{amsart}

\usepackage{times,epsf,amssymb,amsmath,hyperref, rlepsf, color}

\begin{document}

\newtheorem{thm}{Theorem}[section]
\newtheorem{lem}[thm]{Lemma}
\newtheorem{cor}[thm]{Corollary}

\theoremstyle{definition}
\newtheorem{defn}[thm]{\bf{Definition}}

\theoremstyle{remark}
\newtheorem{rmk}[thm]{Remark}

\def\square{\hfill${\vcenter{\vbox{\hrule height.4pt \hbox{\vrule width.4pt height7pt \kern7pt \vrule width.4pt} \hrule height.4pt}}}$}

\newenvironment{pf}{{\it Proof:}\quad}{\square \vskip 12pt}

\title[Embedded $H$-disks]{Embeddedness of the solutions\\ to the $H$-Plateau Problem}
\author{Baris Coskunuzer}
\address{MIT, Mathematics Department, Cambridge, MA 02139}
\address{Boston College, Mathematics Department, Chestnut Hill, MA 02467}
\email{coskunuz@bc.edu}
\thanks{The author is supported by Fulbright Grant, TUBITAK 2219 Grant, and BAGEP award of the Science Academy.}

\maketitle

%% User definitions:

\newcommand{\Si}{S^2_{\infty}({\Bbb H}^3)}
\newcommand{\SI}{S^n_{\infty}({\Bbb H}^{n+1})}
\newcommand{\PI}{\partial_{\infty}}

\newcommand{\BH}{\Bbb H}
\newcommand{\BHH}{{\Bbb H}^3}
\newcommand{\BR}{\Bbb R}
\newcommand{\BC}{\Bbb C}
\newcommand{\BZ}{\Bbb Z}

\newcommand{\e}{\epsilon}

\newcommand{\wh}{\widehat}
\newcommand{\wt}{\widetilde}

\newcommand{\A}{\mathcal{A}}
\newcommand{\C}{\mathcal{C}}
\newcommand{\D}{\mathcal{D}}
\newcommand{\U}{\mathcal{U}}
\newcommand{\V}{\mathrm{E}}
\newcommand{\M}{\mathbf{M}}
\newcommand{\F}{\mathcal{F}}
\newcommand{\I}{\mathcal{I}}
\newcommand{\p}{\mathcal{P}}
\newcommand{\R}{\mathcal{R}}
\newcommand{\B}{\mathbf{B}}
\newcommand{\h}{\mathcal{H}}
\newcommand{\T}{\mathfrak{T}}
\newcommand{\s}{\mathcal{S}}

\begin{abstract}

We generalize Meeks and Yau's embeddedness result on the solutions of the Plateau problem to constant mean curvature disks. We show that any minimizing $H$-disk in an $H_0$-convex domain is embedded for any $H\in[0,H_0)$. In particular, for the unit ball $\B$ in $\BR^3$, this implies that for any $H \in[0,1]$, any Jordan curve in $\partial \B$ bounds an embedded $H$-disk in $\B$.
\end{abstract}

\section{Introduction}

Embeddedness of the solutions to the Plateau problem has always been of great interest in minimal surface theory. However, its natural generalization to  constant mean curvature (CMC) disks has not been addressed yet. While there are many important results on the existence of immersed disks with constant mean curvature or prescribed mean curvature in suitable domains \cite{Hi, HK, Gu, We, BC}, there is virtually no result on the literature on the embeddedness of the solution to this generalized Plateau problem ($H$-Plateau problem).

%Of course, if there is no topological restriction on the surface, the {\em minimizing} $H$-surfaces are always embedded by the regularity theory of Geometric Measure Theory.

In this paper, we will address the existence and the embeddedness of the solutions to the $H$-Plateau problem. In \cite{MY1, MY2, MY3}, Meeks-Yau showed that if $\Gamma$ is a simple closed curve in $\partial \Omega$ where $\Omega$ is a mean convex domain, and if $\Gamma$ is nullhomotopic in $\Omega$, then any area minimizing disk $\D$ in $\Omega$ with $\partial \D=\Gamma$ is embedded. We generalize their result by showing the existence and the embeddedness of $H$-disks in very general domains.

\begin{thm} \label{main1} Let $M$ be an $H_0$-convex manifold with $H_2(M)=\{0\}$. Let $\Gamma$ be a Jordan curve in $\partial M$ such that $\Gamma$ is nullhomotopic in $M$ and separating in $\partial M$. Then, for any $ H\in[0,H_0)$, there exists a minimizing $H$-disk $\D_H$ in $M$ with $\partial \D_H=\Gamma$, and any such minimizing $H$-disk is embedded.
\end{thm}

In particular, for the unit ball $\B$ in $\BR^3$, this result implies that for any $H\in[0,1]$, for any Jordan curve $\Gamma$ in $\partial \B$, the minimizing $H$-disk $\D_H$ with $\partial \D_H=\Gamma$ is embedded. Note that the {\em minimizing} condition for $H$-disks is essential here, e.g. see the examples described in Section \ref{nonexample}.

A corollary of the main result gives a stronger result for the Rellich Conjecture. Rellich conjectured that any Jordan curve $\Gamma$ bounds at least two solutions of $H$-Plateau problem for sufficiently small $H>0$  \cite{BC}. Our result can simply be generalized to give a positive answer to this conjecture for $H_0$-extreme curves. Moreover, these two solutions are both embedded (See Section \ref{Rellich-sec}).

\begin{cor} Let $M$ be an $H_0$-convex manifold with $H_2(M)=\{0\}$. Let $\Gamma$ be a Jordan curve in $\partial M$ such that $\Gamma$ is nullhomotopic in $M$ and separating in $\partial M$. Then, for any $H\in (0,H_0)$, there are two embedded $H$-disk $\D^\pm_H$ in $M$ with $\partial \D^\pm_H=\Gamma$.
\end{cor}

Furthermore in section \ref{Dehn}, we discuss the Dehn's Lemma for $H$-disks, and nonexistence of minimizing $H$-disks for some $H_0$-convex manifolds. In particular, we show the sharpness of Theorem \ref{main1}, by constructing a compact $H_0$ convex $3$-manifold $M$, and a Jordan curve in $\partial M$ which bounds no minimizing $H$-disks in $M$ for some $H\in(0,H_0)$; see Remark \ref{dehn3}.

\vspace{.2cm}

\noindent {\em The outline of proof of the main result:} First, we showed the embeddedness of minimizing $H$-disks for $H_0$-convex domains (Theorem \ref{embed1}) by adapting the techniques in \cite{MY3} with some surgery arguments. The new major tool here is the fact that for $H>0$, the mean curvature vector $\mathbf{H}$ along the minimizing $H$-disk $\Sigma$ always points outside the region $\Omega^-$, where $\Omega^-$ is the region which $\Sigma$ separates from $M$. By using this $H$-concavity property of $\Omega^-$, we can successfully apply the surgery arguments, and show the minimizing $H$-disks must be embedded. As a corollary, we show that the result naturally applies to the small balls in Riemannian manifolds.

Then, to obtain the main result (Theorem \ref{main}), we needed to show the existence of minimizing $H$-disks in $H_0$-convex manifolds. For that, we construct a special "ordered" sequence of embedded disks $\{E_i\}$ in $H_0$-convex manifold $M$ where the limit minimize our variational problem $\I_H$. Then, by using the arguments in \cite{HS}, we restricted our sequence to small balls in $M$, and show that the limit is a smooth $H$-disk in $M$ with the given boundary.

The organization of the paper is as follows: In Section 2, we give basic facts about $H$-disks, and describe the $H$-Plateau problem. Then, in Section 3, we show the embeddedness result for $H$-minimizing disks in $H_0$-convex manifolds. In Section 4, we prove the existence of minimizing $H$-disks, and Theorem \ref{main1}, and show the relevant results mentioned above. Finally in Section 5, we give some concluding remarks.

\subsection{Acknowledgements}

Part of this research was carried out at Massachusetts Institute of Technology during my visit. I would like to thank them for their great hospitality.

\section{Preliminaries}

In this section, we will overview the basic results and notions which we use in the following sections. First, we recall the CMC disks, and the related variational problem. We will be working with the immersions from the unit disk to a Riemannian $3$-manifold $M$. In particular,  when we use the term immersion, we mean a map $u:D^2\to M$ such that $u\in \C^0(\overline{D^2})\cap \C^2(D^2)$. Let $\Sigma_u=u(D^2)$ be an immersed disk in $M$. If the mean curvature at $u(p)\in \Sigma_u$ is equal to $H$ for any $p\in int(D^2)$, we will call $\Sigma_u$ an immersed $H$-disk in $M$.

\subsection{$H$-Plateau Problem} \label{H-Plateau} \

In the references mentioned in the introduction \cite{Hi, HK, Gu, We, BC}, the existence of an $H$-disk bounding a given simple closed curve in $\BR^3$ (or $M$) was studied extensively. Furthermore, to show the existence of an $H$-disk, they considered the following variational problem.

For simplicity, we will take the unit ball $\B$ as the ambient space. Let $\Gamma$ be a simple closed curve in $\B$. Define the space of immersions $$X_\Gamma=\{u:D^2\to \B \ | \ u\in \C^0(\overline{D^2})\cap \C^2(D^2) \mbox{ and } u(\partial D^2)=\Gamma\}$$

Define the operator $F_H:X_\Gamma\to\BR$ as follows: $$\F_H(u)=\int_{D^2} |u_x|^2+ |u_y|^2 +\dfrac{4}{3}H [u\cdot(u_x\times u_y)] \ dxdy$$

In other words, $\F_H(u)=E(u)+\dfrac{4}{3}H W(u)$ where $$E(u)=\int_{D^2} |u_x|^2+|u_y|^2 \ dxdy \ \mbox{ and } \ W(u)=\int_{D^2} [u\cdot(u_x\times u_y)] dxdy$$
The critical points of $\F_H$ are conformal, $H$-harmonic maps. Hence, for such a $u$, if $\Sigma_u=u(D^2)$ is the image, then $\Sigma_u$ is an immersed $H$-disk. They showed the existence of an immersed $H$-disk $\Sigma_u$ with $\partial \Sigma_u =\Gamma$ by proving the existence of an immersion $u$ which {\em minimizes} $F_H$.

In particular, when $H=0$, this is the classical Plateau problem. Hence, minimizing $F_H$ in $X_\Gamma$ is considered as a natural generalization of the classical Plateau problem. Hence, we will call this problem as {\em $H$-Plateau problem}. Now, we rephrase the definition above.

\begin{defn} {\bf i.} Let $u:D^2\to M$ be an immersion, and let $H\geq 0$. If $\Sigma=u(D^2)$ is a critical point of $\F_H$ for any variation, we call $\Sigma$ an \textit{$H$-disk}. Equivalently, $\Sigma$ has constant mean curvature $H$ at every point.

{\bf ii.} An immersed $H$-disk $\Sigma$ is a \textit{minimizing $H$-disk} if $\Sigma$ is the absolute minimum of the functional $\F_H$ among the immersed disks with the same boundary.
\end{defn}

\noindent \textbf{Notation:} From now on, we will call CMC surfaces with mean curvature $H$ as \textit{$H$-surfaces}. All the surfaces are assumed to be orientable unless otherwise stated.

\subsection{A Reformulation of the $H$-Plateau Problem} \

There is another formulation of the $H$-Plateau problem in a slightly different way for the extreme curves. To keep the paper self contained, we recall the relation between these formulations. Let $\Gamma$ be a Jordan curve in $\partial \B$ and let $\partial \B-\Gamma= E^+\cup E^-$. Let $\Omega^-_u$ be the open region in $\B-\Sigma_u$ which $\Sigma_u\cup E^-$ separates from $\B$, i.e. $\partial \overline{\Omega_u^-}=\Sigma_0\cup E^-$.

For a given $u\in \C^0(\overline{D^2})\cap \C^2(D^2)$, consider the variational problem: $$\I_H(u)=Area(\Sigma_u)+2H.Vol(\Omega^-_u)$$

Here, $Vol(.)$ represents the \underline{oriented} volume of the region.

\begin{lem} \label{variation} Let $\F_H$ ,$\I_H$, $u$ and $\Sigma_u$ be as described above. If $u_0$ minimizes $\F_H$, then $\Sigma_{u_0}$ minimizes $\I_H$. Conversely, if $u_0$ is conformal, and $\Sigma_{u_0}$ minimizes $\I_H$, then $u_0$ minimizes $\F_H$.
\end{lem}

\begin{pf}  Let $u_0$ minimize $\F_H$. For any $u$, $2 Area(\Sigma_u)\leq E(u)$ with the equality if and only if $u$ is conformal. Since the minimizer $u_0$ for $\F_H$ is conformal by \cite{Hi}, we have $2 Area(\Sigma_{u_0})= E(u_0)$

Now, notice that minimizing $W(u)$ would be the same with minimizing $Vol(\Omega^-_u)$ by the divergence theorem. In particular,

{\small $$3Vol(\Omega^-_u)= \int_{\Omega^-_u} \nabla \cdot \mathbf{x} \ dV= \int_{\partial \Omega^-_u} \mathbf{x} \cdot \mathbf{N} \ dS= \int_{\Sigma_u}\mathbf{x}\cdot\mathbf{N} \ dS + \int_{E^-}\mathbf{x}\cdot\mathbf{N} \ dS=  W(u)+ C_0$$}

where $C_0$ is the constant corresponding the quantity $\int_{E^-}\mathbf{x}\cdot\mathbf{N} dS$. Hence, if $u$ minimizes $\F_H$, then

\begin{align*}
\F_H(u) ={}&  E(u)+\dfrac{4}{3}H.W(u)=2Area(u)+4H(Vol(\Omega^-_u)- C_1){}\\
    ={} &  2[Area(u)+2H.(Vol(\Omega^-_u)-C_1)]=2\I_H(u)-C_2
\end{align*}

%$\F_H(u)=E(u)+\dfrac{4}{3}H.W(u)=E(u)+4H[Vol(\Omega^-_u)- C_1]= 2[Area(u)+2H.(Vol(\Omega^-_u)-C_1)]=2\I_H(u)-C_2$

where $C_i$ are the corresponding constants. By the existence of isothermal coordinates, this implies if $u$ minimizes $\F_H$, then it also minimizes $\I_H$.

Conversely, for conformal $u$, if $\Sigma_u$ minimizes $\I_H$, then $u$ minimizes $\F_H$ by the identity above.
\end{pf}

Throughout the paper, we will use the variational problem described by $\I_H$. The advantage of working with $\I_H$ is that it is independent of the parametrization of $\Sigma_u$. Note that these equivalent variational characterizations $\I_H$ and $\F_H$ for $H$-surfaces differs at second variations \cite{Gr}.

\begin{rmk} [Minimizing $H$-surfaces] Note that the variational problem $\I_H$ is widely used to show existence of $H$-surfaces (with no topological restriction) with geometric measure theory techniques. In particular, the compactness and interior regularity results of geometric measure theory gives {\em smoothly embedded} $H$-surfaces with given boundary curves by minimizing $\I_H$ in a suitable setting. However, in this approach, there is no control on the topology of the $H$-surface.

If we restrict the topology of the surface to the disk, for a given general simple closed curve $\gamma$ in $M$, there might be self intersections of the minimizing $H$-disks (if exists) bounding $\gamma$. Even in the area minimizing case ($H=0$), it is easy to construct many simple closed curves in any ambient space where the area minimizing disk is not embedded. In this paper, we will address to this question.
\end{rmk}

\vspace{.3cm}

\subsection{$H_0$-convex domains} \

These are mean convex domains where $H_0$ denotes the lower bound for the mean curvature of the boundary.

\begin{defn} [$H_0$-convex domains] Let $\Omega$ be a compact $3$-manifold with piecewise smooth boundary. We call $\Omega$ an {\em $H_0$-convex domain} if

\begin{itemize}
\item The mean curvature vector $\mathbf{H}$ always points towards inside $\Omega$ along the smooth parts of $\partial \Omega$,

\item  The mean curvature $|\mathbf{H}(p)|\geq H_0$ for any smooth point  $p\in\partial \Omega$,

\item Along the nonsmooth parts of $\partial \Omega$, the inner dihedral angle is less than $\pi$.

\end{itemize}
\end{defn}

%Let $M$ be a compact Riemannian $3$-manifold with boundary, and let $H_0\geq 0$. Then, $M$ is called {\em $H_0$-convex} if for any point $p\in\partial M$, the mean curvature vector $\mathbf{H}(p)$ always points inside $M$, and the mean curvature $|\mathbf{H}(p)|\geq H_0$.

This is a natural generalization of {\em sufficiently convex} or {\em mean convex} notion used in \cite{MY2}, \cite[Section 6]{HS}. For further generalization, we will call curves in the boundary of a $H_0$-convex domains as {\em $H_0$-extreme curves}.

\vspace{.2cm}

\noindent {\em Examples:} Any mean convex manifold is automatically $0$-convex with respect to this definition. Similarly, an $R$-ball in $\BR^3$ is $\frac{1}{R}$-convex.

\begin{rmk} \label{homreg} [Noncompact $H_0$-convex domains] The compactness assumption on the $H_0$-convex domain is not really essential here. One can replace the compactness assumption with homogeneously regular assumption for noncompact manifolds. Here, $M$ is {\em homogeneously regular} means that the sectional curvatures are bounded above, and injectivity radii are bounded below in $M$ \cite{HS}. Hence, we will call a noncompact $3$-manifold $M$ with boundary as $H_0$-convex domain if $M$ is homogeneously regular, and $\partial M$ is $H_0$-convex (satisfying the conditions above). All the techniques and results in this paper naturally apply to this noncompact case \cite[Section 6]{HS}.
\end{rmk}

\begin{rmk} \label{0-convex} With this notation, $0$-convex domains exactly correspond to the mean convex domains in \cite{MY2}. Hence, Meeks-Yau's embeddedness result can be stated as follows: If $\Gamma$ is a simple closed curve in $\partial \Omega$ where $\Omega$ is a $0$-convex domain, and if $\Gamma$ is nullhomotopic in $\Omega$, then any minimizing $0$-disk $\Sigma$ in $\Omega$ with $\partial \Sigma=\Gamma$ is embedded. In this paper, we will generalize this result by replacing $0$ in the statement above with $H$ for $H_0$-convex domains where $ H\in[0,H_0)$ (Theorem \ref{main}).
\end{rmk}

%If $\Gamma$ is a simple closed curve in $\partial \Omega$ where $\Omega$ is a $H_0$-convex domain, and if $\Gamma$ is nullhomotopic in $\Omega$, then any minimizing $H$-disk $\Sigma$ in $\Omega$ with $\partial \Sigma=\Gamma$ is embedded

\subsection{Existence of Solutions to the $H$-Plateau Problem} \

There are many results on the existence of immersed $H$-disks in $H_0$-convex domains in the literature \cite{Hi, HK, Gu, We}. These existence results were obtained mainly by solving the $H$-Plateau problem $\F_H$ described in Section \ref{H-Plateau}. The basic case is the round balls $\B_R$ in $\BR^3$.

\begin{lem} \label{H-disk1} \cite{Hi} Let $\B_R$ is a closed ball in $\BR^3$ of radius $R$. Let $\Gamma$ be a Jordan curve in $\B_R$. Then, for any $H\in [0,\frac{1}{R}]$, there exists an immersion $u_H:D^2\to \B_R$ such that $u_H(\partial D^2)=\Gamma$ and the image $u_H(D^2)$ is a minimizing $H$-disk.
\end{lem}

This result is sharp, as for the round equator circle in $\partial \B_R$, there is no solution to the $H$-Plateau problem in $\B_R$ for any $H>\frac{1}{R}$.

A generalization of this result to the Riemannian manifolds came a few years later.

\begin{lem} \label{H-disk2} \cite{Gu, HK} Let $M$ be a homogeneously regular $3$-manifold. For any $H_0>0$, for any point $p\in M-\partial M$, there exists a sufficiently small  $\e>0$ such that for any simple closed curve $\Gamma$ in $\overline{B_\e(p)}$, there exists an immersed minimizing $H$-disk $\Sigma_H$ in $\overline{B_\e(p)}$ with $\partial \Sigma_H = \Gamma$ where $ H\in [0,H_0]$.
\end{lem}

The statements of this result in \cite{Gu}, and \cite{HK} are more general, but for our purposes, we will just take this simpler version here. Indeed, for any point $p\in M-\partial M$, one can find an immersed $H$-disk bounding any given curve in $B_r(p)$ for $H\in [0,\sqrt{K_0}\cot(\sqrt{K_0}r)]$ where $K_0$ is an upper bound for the sectional curvature of the ambient manifold $M$ \cite{Gu}.

We finish this section with the following lemma, which is known as the maximum principle.

\begin{lem} \cite{Gu} \label{maximum} (Maximum Principle) Let $\Sigma_1$ and $\Sigma_2$ be two surfaces in a Riemannian $3$-manifold which intersect at a common point tangentially. Let $H_i$ be the (signed) mean curvature of $\Sigma_i$ at the common point with respect to the same normal vector $\mathbf{N}$, i.e. $\mathbf{H}_i=H_i\mathbf{N}$. If $\Sigma_2$ lie in positive side (the normal vector $\mathbf{N}$ direction) of $\Sigma_1$ nearby the common point, then $H_1$ is strictly less than $H_2$, i.e. $H_1 < H_2$.
\end{lem}
%Let $\Sigma_1$ and $\Sigma_2$ be two surfaces in a Riemannian $3$-manifold which intersect at a common point tangentially. If $\Sigma_2$ lies in positive side (mean curvature vector direction) of $\Sigma_1$ around the common point, then $H_1$ is strictly less than $H_2$ ($H_1 < H_2$) where $H_i$ is the mean curvature of $\Sigma_i$ at the common point.

\section{Embeddedness of minimizing $H$-disks in $H_0$-convex domains}

In this section, we will show that the {\em minimizing $H$-disks} in $H_0$-convex domains are embedded by generalizing the ideas in \cite{MY3} to CMC setting.

First, we need a lemma proving the structure of the self intersection set of immersed $H$-disks. Let $M$ be an $H_0$-convex domain, and let $\Gamma$ be a simple closed curve in $\partial M$. Define the space of immersions from a disk to $M$ as follows:
$$X_\Gamma=\{u:D^2\to M \ | \ u\in \C^0(\overline{D^2})\cap \C^2(D^2) \mbox{ and }  u|_{\partial D^2} \mbox{ is an embedding onto } \Gamma\}$$

For a given immersion $u\in X_\Gamma$, define the self-intersection set $Y_u$ as follows: $Y_u=\{p\in \overline{D^2} \ | \ \exists q\neq p \mbox{ with } u(p)=u(q)\}$.

\begin{lem} \cite{MY3} \label{self} Let $M$ be an $H_0$-convex domain, and let $\Gamma$ be a simple closed curve in $\partial M$. Fix $ H \in[0,H_0)$. Let $u\in X_\Gamma$ be a conformal, $H$-harmonic map. Then, $Y_u$ is a $1$-complex in $int(D^2)$ where every vertex is joined by at least two edges.
\end{lem}

\begin{pf} The proof follows from \cite[Lemma 3]{MY3}, and its corollary. Meeks and Yau showed this result for conformal harmonic immersions, but the same proof works exactly by using the maximum principle (Lemma \ref{maximum}) for $H$-harmonic maps instead of the maximum principle for minimal surfaces. In particular, as $M$ is $H_0$-convex, the $H$-disk $\Sigma_u=u(D^2)$ cannot intersect $\partial M$ in the interior by Lemma \ref{maximum}, i.e. $\Sigma_u\cap \partial M=\Gamma$. This shows that $Y_u\subset int(D^2)$. Again by the maximum principle, there is no isolated points in $Y_u$. Similarly, there cannot be a curve with an endpoint in $Y_u$. Finally, again by the maximum principle, for any two disjoint open subsets $O_1$ and $O_2$ in $D^2$, $u(O_1)\neq u(O_2)\subset \Sigma_u$. Hence, $Y_u$ is a $1$-complex where every vertex is joined by at least two edges.
\end{pf}

Informally, the above lemma says that $Y_u=\gamma_1\cup\gamma_2\cup...\cup\gamma_m$ where $\gamma_i$ is a closed curve in $int(D^2)$, and $\gamma_i\cap\gamma_j$ is either empty or a finite set of points (vertices).

Now, we show that any immersed solution to $H$-Plateau problem in an $H_0$-convex manifold must be embedded.

%Now, we restrict ourselves to the unit ball $\B$. Later in this section, we give a generalization of the following result for small balls in Riemannian $3$-manifolds.

\begin{thm} \label{embed1} Let $M$ be a compact (or homogenously regular) $H_0$-convex manifold with $H_2(M)=\{0\}$. Let $\Gamma$ be a Jordan curve in $\partial M$ such that $\Gamma$ is nullhomotopic in $M$ and separating in $\partial M$. Fix $H\in [0,H_0)$. If $\Sigma$ is a minimizing $H$-disk in $M$ with $\partial \Sigma=\Gamma$, then $\Sigma$ is embedded.
\end{thm}

\begin{pf} We will use the notation set in section \ref{H-Plateau}. For $H=0$, the result follows from \cite{MY2}. So, we assume $H\in(0,H_0)$.

Let $u_0:D^2\to M$ be with $u_0\in \C^0(\overline{D^2})\cap \C^2(D^2)$ and $u_0(\partial D^2)=\Gamma$ such that $u_0$ minimizes $\F_H$. In particular, let the space of immersions $X_\Gamma:\{u:D^2\to \B \ | \ u\in \C^0(\overline{D^2})\cap \C^2(D^2) \mbox{ and } u(\partial D^2)=\Gamma\}$ be as above. Then $\F_H(u_0)=\inf_{X_\Gamma} \{\F_H(u)\}$. Hence, $u_0$ is conformal, and $\Sigma_0=u_0(\overline{D^2})$ is an immersed $H$-disk in $M$ with $\partial \Sigma_0 = \Gamma$. 

Let $S$ be the component in $\partial M$ with $\Gamma\subset S$. Since $\Gamma$ is separating in $\partial M$, let $S-\Gamma=int(S^-)\cup int(S^+)$ with $\partial S^\pm=\Gamma$. Then, $\Sigma_0\cup S^-$ would be a closed surface in $M$. Since $H_2(M)$ is trivial, $\Sigma_0\cup S^-$ bounds a domain $\Omega_0^-$ in $M$. Then by using $\Omega_0^-$, we can define $I_H(\Sigma_0)=|\Sigma_0|+2H\|\Omega_0^-\|$ where $|.|$ represents the area, and $\|.\|$ represents the volume. as before. Then, by Lemma \ref{variation}, $\Sigma_0$ minimizes $\I_H$, i.e. $\I_H(\Sigma_0)=\inf_{X_\Gamma} \I_H(\Sigma_u)$.

%We claim that $u_0$ is an embedding, and $\Sigma_0$ is an embedded $H$-disk (no self intersections). Assume that $\Sigma_0$ is not an embedding. By Lemma \ref{maximum}, $\Sigma_0\cap\partial \B=\Gamma$ as $H\leq 1$. Furthermore, again by the maximum principle, $\Sigma_0$ has no isolated points of self intersections.

%Let $\Omega_0^-$ be the open region in $\B$ which $\Sigma_0\cup E^-$ separates from $\B$, i.e. $\partial \overline{\Omega_0^-}=\Sigma_0\cup E^-$. By the assumption, $\I_H(\Sigma_0)=|\Sigma_0|+2H \|\Omega_0^-\|=\inf_{X_\Gamma} \I_H(\Sigma_u)$ where $|.|$ represents the area, and $\|.\|$ represents the volume.

\vspace{.2cm}

\noindent {\bf Step 1:} [H-concavity of $\Omega_0^-$ ] For any smooth point $p\in \Sigma_0-\Gamma$, the mean curvature vector $\mathbf{H}(p)$ points outside $\Omega_0^-$.

\vspace{.2cm}

\noindent {\em Proof of the Step 1:} Assume otherwise, and let $q$ be a point $\Sigma_0$ where the mean curvature vector points inside $\Omega_0^-$. Let $\e>0$ be sufficiently small so that $D_q=\B_\e(q)\cap \Sigma_0$ is a smooth disk. Then $\gamma_q=\partial D_q$ is a simple closed curve in $\partial \B_\e(q)$. Let $\Omega_q=\B_\e(q)\cap \overline{\Omega_0^-}$. By assumption, the mean curvature vector $\mathbf{H}$ points inside $\Omega_q$ along $D_q$. Hence, $\Omega_q$ is a mean convex domain. Then, there exists an area minimizing disk $D'_q$ in $\Omega_q$ with $\partial D'_q =\gamma_q$ by \cite{MY2}. Note that $|D_q'|<|D_q|$.

Let $\Delta_q$ be the region in $\Omega_q$ with $\partial \Delta_q= D_q\cup D_q'$. Let $\Sigma_0'=(\Sigma_0-D_q)\cup D'_q$, and $\Omega_0'=\Omega_0^--\Delta_q$. Clearly, $|\Sigma_0'|<|\Sigma_0|$ and $\|\Omega_0'\|< \|\Omega_0^-\|$. Then, $\I_H(\Sigma_0')=|\Sigma_0'|+2H \|\Omega_0'\|<|\Sigma_0|+2H \|\Omega_0^-\|=\I_H(\Sigma_0)$. By smoothing out $\Sigma_0'$ along $\gamma_q$, we can get $\Sigma_0''$ in $X_\Gamma$ with $\I_H(\Sigma_0'')<\I_H(\Sigma_0)$. This is a contradiction as $\I_H(\Sigma_0)=\inf_{X_\Gamma} \I_H(\Sigma_u)$. Hence, Step 1 follows. \hfill $\Box$

\vspace{.2cm}

Note that even though $\Sigma$ has self-intersection, $\Omega^-$ can still be $H$-concave as illustrated in Figure \ref{h-concave}-right. In particular, in Figure \ref{h-concave}, the ambient space is the unit ball $\B$ in $\BR^3$. In this 2D figure, blue arcs represents the the immersed $H$-disks $\Sigma_i$, and the shaded region represents the domains $\Omega_i^-$. The red dots in the boundary represents the Jordan curves $\Gamma_i$, and the red dots in the middle represents the self intersection curves. Hence, $\Gamma_i$ is a Jordan curve in $\partial \B$, and  $\Sigma_i$ is an the immersed disk in $\B$ with $\partial \Sigma_i=\Gamma_i$.

In the left $\Omega_1^-$ is not $H$-concave, as at some points (along $\partial X$) on $\Sigma_1$, $\mathbf{H}$ points inside $\Omega_1^-$. Hence, by Step 1, there cannot be a minimizing $H$-disk with self-intersection as in Figure \ref{h-concave}-left. In the right, $\Omega_2^-$ represents a region with two components, where one of them is a solid torus, and the other one is a ball. Notice that even though $\Sigma_2$ is not embedded, $\Omega_2^-$ is still $H$-concave. In particular, $H$-concavity of $\Omega^-$ is a very powerful tool for $H>0$ case, and it rules out self-intersections with even number of preimages (See Remark \ref{h-concavity}).

\begin{figure}[t]
\begin{center}
$\begin{array}{c@{\hspace{.2in}}c}

\relabelbox  {\epsfysize=2in \epsfbox{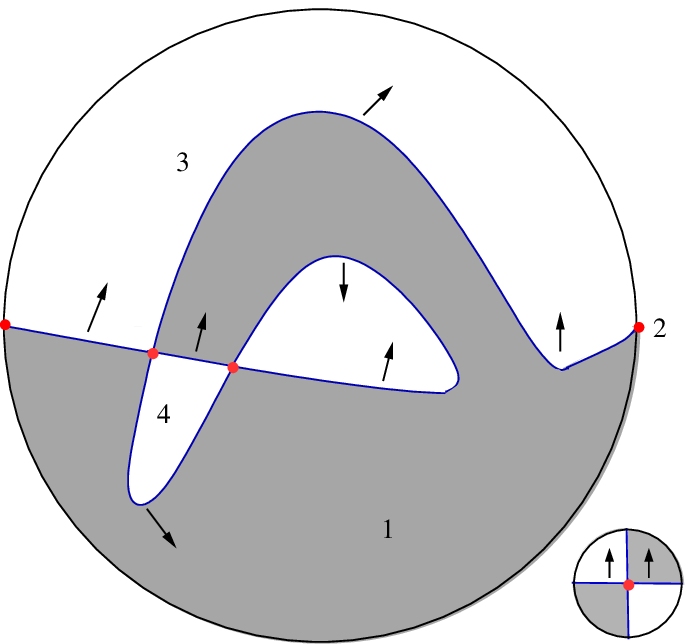}} \relabel{1}{\tiny $\Omega_1^-$} \relabel{2}{\tiny $\Gamma_1$} \relabel{3}{\tiny $\Sigma_1$} \relabel{4}{\tiny $X$} \endrelabelbox &

\relabelbox  {\epsfysize=2in \epsfbox{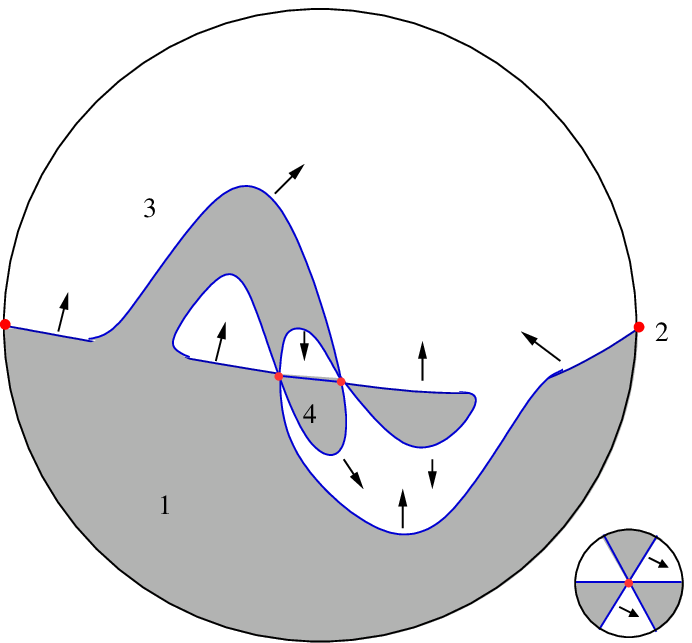}} \relabel{1}{\tiny $\Omega_2^-$} \relabel{2}{\tiny $\Gamma_2$} \relabel{3}{\tiny $\Sigma_2$} \relabel{4}{\tiny $X$} \endrelabelbox \\
\end{array}$

\end{center}
\caption{\label{h-concave} \footnotesize The arrows show the mean curvature vector $\mathbf{H}$ along the immersed disk $\Sigma_i$. In the left, $\Omega_1^-$ is not an $H$-concave region. In the right, $\Omega_2^-$ is an $H$-concave region. For further details, see Remark \ref{h-concavity}.}
\end{figure}

%Now, define set of self-intersection points $Y=\{p\in\Sigma_0 \ | \ |u_0^{-1}(p)|\geq 2 \}$ in $\Sigma_0$. Recall that $u_0:D^2\to \B$ be conformal $H$-harmonic immersion, and hence $\Sigma_0=u_0(D^2)$ is an immersed $H$-disk in $\B$ with $\partial \Sigma_0=\Gamma$. Furthermore, since $\Sigma_0$ is an $H$-disk with $0\leq H\leq 1$, $u_0(int(D^2))\cap\partial \B=\emptyset$ by maximum principle. This implies that $Y\cap \Gamma=\emptyset$.

%Let $Z=u_0^{-1}(Y) \subset D^2$. Notice that $Z$ is a closed subset in $D^2$, and $Z\cap \partial D^2=\emptyset$. Furthermore, there is no isolated point in $Z$ by the maximum principle. Finally, there cannot be a closed arc in $Z$ as $Z\cap \partial D^2=\emptyset$, and $u_0$ is an immersion. This shows that $Z$ is a union of smooth simple closed curves, and hence $Y$ is a union of smooth simple closed curves, i.e. $Y=\wh{\gamma}_1\cup\wh{\gamma}_2\cup ... \cup \wh{\gamma}_n$. Furthermore, for any fixed $i$, $u_0^{-1}(\wh{\gamma_i})=\bigcup_{j=1}^{m_i}\gamma_i^j\cup X_i$ which consists of $m_i\geq 2$ pairwise disjoint simple closed curves $\{\gamma_i^j\}$ in $Z\subset D^2$ and $X_i$ is a finite set of points which belongs to other curves $\gamma_{i'}^{j'}$. Again by the maximum principle, $\gamma_i^j\cap\gamma_{i'}^{j'}$ is either empty or a finite set of points, i.e. two curves in $Z$ cannot coincide along an arc.

Let $Y_0$ be the self-intersection set of $u_0$. Recall that $Y_0$ is a collection of closed curves in $int(D^2)$, possible intersecting each other at a finite set of points (vertices).

\vspace{.2cm}

\noindent {\bf Step 2:} [Surgery] $\Sigma_0$ is an embedding.

\vspace{.2cm}

\noindent {\em Proof of the Step 2:} Assume that the self-intersection set $Y_0\neq \emptyset$. Then, by the tower construction of Papakyriakopoulous \cite[Theorem 4]{MY3}, for a lift of $N(\Sigma_0)$, there are two simple closed curves $\gamma^+$ and $\gamma^-$ in $Y_0$ where $u(\gamma^+)=u(\gamma^-)=\wh{\gamma}\subset \Sigma_0$, and $\gamma^\pm=\partial D^\pm$ with $Y_0\cap int(D^\pm)=\emptyset$.

Here, $D^+$ and $D^-$ are the closed subdisks in $D^2$ with $\partial D^\pm =\gamma^\pm$. Let $u(D^\pm)=\wh{D^\pm}\subset \Sigma_0$ where $\partial \wh{D^\pm}=\wh{\gamma}$. Now, we define another map $w:D^2\to M$ by swaping the disks $\wh{D^+}$ and $\wh{D^-}$ in $\Sigma_0$. In particular, $w=u_0$ on $D^2-int(D^+\cup D^-)$. Let $\varphi:D^+\to D^-$ be a diffeomorphism. Then, define $w=u\circ \varphi$ on $D^+$. Similarly, let $w=u\circ \varphi^{-1}$ on $D^-$. Hence, $w$ is an immersion on $D^2$ except along $\gamma^+\cup\gamma^-$. Let $\Sigma'=w(D^2)$.

Clearly, $\Sigma_0=\Sigma'$ by construction. Similarly, if $\Omega'$ is the region in $\B$ with $\partial \Omega'= \Sigma'\cup E^-$, then $\Omega'=\Omega_0$. By the $H$-concavity of $\Omega_0$ (Step 1), when we swap the $\wh{D^+}$ and $\wh{D^-}$ in $\Sigma_0$, the normal vectors in both disks still point outside $\Omega_0$. Hence, the oriented volume of $\Omega'$ and $\Omega_0$ are the same. This implies $\I_H(\Sigma_0)=\I_H(\Sigma')$.

However, the new map $w$ is not smooth, and we have a folding curve $\wh{\gamma}$ in $\Sigma'$. Hence, by \cite[Lemma 7]{MY3}, if we push $\Sigma'$ towards $\Omega'$ along $\wh{\gamma}$, we decrease the area of $\Sigma'$, and we decrease the volume of $\Omega'$ at the same time. Hence, we can smooth out $\Sigma'$ by pushing along $\wh{\gamma}$, and get an immersion $\wt{\Sigma}=\wt{w}(D^2)$ where $\wt{w}\in X_\Gamma$. However, $I_H(\wt{\Sigma})<\I_H(\Sigma_0)$ while $\wt{w}\in X_\Gamma$. This contradicts to the fact that $\I_H(\Sigma_0)$ is the infimum in $X_\Gamma$.
\end{pf}

\begin{rmk} \label{h-concavity} [$H$-concavity of $\Omega^-$] Notice that Step 1 in the proof above implies any arc in the self intersection set $Y_0$ cannot have even number of preimages in $D^2$ under $u_0$. In particular, the Figure \ref{h-concave}-left shows that if the curve $\gamma\subset Y_0$ (red dots in the boundary of the region $X$) has two preimages, then the mean curvature  vector $\mathbf{H}$ along $\Sigma_1$ must point inside $\Omega_1^-$ in one side of $\gamma$ which contradicts to the $H$-concavity of $\Omega^-$. As a result, in the left, $\mathbf{H}$ points inside $\Omega_1^-$ along $\partial X$. See the little disk in the left picturing the neighborhood of $\gamma$ to see why $\Omega_1^-$ cannot be $H$-concave.

Hence, any arc $\alpha$ in $Y_0$ must have odd number of preimages under $u_0$ to keep $\Omega^-$ $H$-concave. In the Figure \ref{h-concave}-right, the curve $\alpha\subset Y_0$ (red dots in the boundary of $X$) have 3 preimages. Then, the mean curvature  vector $\mathbf{H}$ along $\Sigma_2$ points outside $\Omega_2^-$ in both sides of $\alpha$ and along $\partial X$, which keeps $\Omega_2^-$ $H$-concave. See the little disk in the right picturing the neighborhood of $\alpha$. In Step 2, we rule out this case, and show the embeddedness of the minimizing $H$-disk.

In particular, $H$-concavity of $\Omega^-$ gives a strong tool to show embeddedness of the minimizing $H$-disks. This is the main difference with the area minimizing ($H=0$) case, and it gives a great advantage for $H>0$ case.
\end{rmk}

The following corollaries are straightforward.

\begin{cor} \label{embed-ball} Let $\B_R$ be a closed ball of radius $R$ in $\BR^3$. Let $\Gamma$ be a Jordan curve in $\partial \B_R$, and let $H \in[0,\frac{1}{R}]$. Then, there exists a minimizing $H$-disk $\Sigma_H$  in $\B_R$ with $\partial \Sigma_H=\Gamma$, and any such minimizing $H$-disk is embedded.
\end{cor}

\begin{pf} Existence follows from Lemma \ref{H-disk1}. Clearly, $\B_R$ satisfies the conditions of Theorem \ref{embed1} for an $\frac{1}{R}$-convex domain with $H_2(\B_R)=\{0\}$. Furthermore, $\partial \B_R$ is a sphere where any Jordan curve $\Gamma$ is separating. The proof follows by Theorem \ref{embed1}.
\end{pf}

\begin{cor} \label{embed2} Let $M$ be a homogeneously regular $3$-manifold. For any $H_0>0$, for any point $p\in M-\partial M$, there exists a sufficiently small  $\e>0$ such that for any simple closed curve $\Gamma$ in $\partial B_\e(p)$  and $H \in [0,H_0)$, there exists a minimizing $H$-disk $\Sigma_H$  in $B_\e(p)$ with $\partial \Sigma_H=\Gamma$, and any such minimizing $H$-disk is embedded.
\end{cor}

\begin{pf} Existence of the minimizer follows from Lemma \ref{H-disk2}. For sufficiently small $\e>0$, $\B_\e$ is an $H_0$-convex domain with $H_2(\B_\e)=\{0\}$. Furthermore, $\partial \B_\e$ is a sphere where any Jordan curve is separating. The proof follows by Theorem \ref{embed1}.
	
\end{pf}

\section{Existence of minimizing $H$-disks in $H_0$-convex manifolds}

%\vspace{-.2cm}

Now, we prove the main result of the paper, a generalization of Meeks and Yau's celebrated embeddedness result to CMC case (See Remark \ref{0-convex}).

In Theorem \ref{embed1}, we proved that any minimizing $H$-disks must be embedded in $H_0$-convex domains. In this section, we prove the existence of minimizing $H$-disks in these domains, and obtain the main result.

%We first show the existence of $H$-disks in a very general setting, and then generalize Meeks and Yau's celebrated embeddedness result to CMC case (See Remark \ref{0-convex}).

%In their celebrated papers \cite{MY1,MY2,MY3}, Meeks and Yau  showed the embeddedness of the area minimizing disks in mean convex manifolds (See Remark \ref{0-convex}). We generalize their result to CMC case as follows.

%In particular, this is considered as {\em Dehn's Lemma} in low dimensional topology terms.

\begin{thm} \label{main} Let $M$ be a compact (or homogenously regular) $H_0$-convex manifold with $H_2(M)=\{0\}$. Let $\Gamma$ be a Jordan curve in $\partial M$ such that $\Gamma$ is nullhomotopic in $M$ and separating in $\partial M$.  Then, for any $ H\in[0,H_0)$, there exists a minimizing $H$-disk $\Sigma_H$ in $M$ with $\partial \Sigma_H=\Gamma$, and any such minimizing $H$-disk is embedded.
\end{thm}

\begin{pf} Notice that embeddedness of such minimizing $H$-disks in $M$ follows from Theorem \ref{embed1}. Hence, we only need to show the existence of such a minimizer to finish the proof.
	
Fix $H\in [0,H_0)$. Define the space of embedded disks in $M$ bounding $\Gamma$, $Y_\Gamma=\{ D \subset M \ | \ D \mbox{ is an embedded disk with } \partial D=\Gamma \}$. Similarly, define the space of immersed disks in $M$ bounding $\Gamma$, $X_\Gamma=\{ D \subset M \ | \ D \mbox{ is an immersed }$  $\mbox{disk with } \partial D=\Gamma \}$.

Let $S$ be the component of $\partial M$ containing $\Gamma$. Since $\Gamma$ is separating in $\partial M$, $\Gamma$ separates $S$ into two parts, say $S^\pm\subset S$ with $\partial S^\pm=\Gamma$ and $S-\Gamma=int(S^+)\cup int(S^-)$. For any embedded disk $D\in Y_\Gamma$, $D\cup S^-$ is a closed surface in $M$. Since $H_2(M)=\{0\}$, it separates a domain $\Omega^-$ in $M$ with $\partial \overline{\Omega^-} = D\cup S^-$. Now, define the map $\I_H: Y_\Gamma \to \BR$ such that $\I_H(D)=Area(D)+ 2HVol(\Omega^-)$. 
Let $b_\Gamma=\inf \{\I_H(D) \ | \ D\in Y_\Gamma\}$. Notice that $b_\Gamma=\inf_{Y_\Gamma} \{\I_H(D)\}=\inf_{X_\Gamma} \{\I_H(D)\}$ by Theorem \ref{embed1}.

\vspace{.2cm}

\noindent {\em Outline of the Proof:} Our aim is to construct a special ordered sequence of disks $\{E_i\}$ with $\I_H(E_i)\searrow b_\Gamma$. To define such sequence, we first construct a lower barrier $N$, and upper barrier $\wh{\Omega}$ in $M$ where $S^-\subset N\subset \wh{\Omega}$ (Step 1 and 2). If $E$ is an embedded disk in $Y_\Gamma$, let $\Delta$ be the component of $M-E$ with $\partial \Delta\supset int (S^-)$. Then, we show that there exists a sequence of pairwise disjoint disks  $\{E_i\}$ such that $N\subset \Delta_i\subset \Delta_j\subset \wh{\Omega}$ for $i<j$ and $\I_H(E_i)\searrow b_\Gamma$ (Step 3). Then, we show that the limit of this sequence is indeed a smooth $H$-surface $\Sigma$ in $M$ with $\partial \Sigma = \Gamma$ (Step 4). To show this, we restrict ourselves to a small ball $B_\e(p)$ around a limit point $p\in \Sigma$, and we modify the sequence $\{E_i\}$ by replacing $E_i\cap B_\e(p)$ with a minimizing $H$-disk $U_i$ given by Corollary \ref{embed2}. By following \cite{HS} and ordered structure of the sequence $\{E_i\}$, we show that the limit of minimizing $H$-disks $U_i$ is a smoothly embedded $H$-disk $U$. By using the property $\I_H(E_i)\searrow b_\Gamma$, we conclude that $U= \Sigma\cap B_\e(p)$ which shows that $\Sigma$ is a smooth $H$-surface. Finally, in Step 5, we prove that $\Sigma$ is indeed a $H$-disk.

%First, we will construct an open neighborhood $N$ of $S^-$ in $M$ which will act as a barrier near $S^-$. Then, we will define a special sequence of embedded disks $\{D_i\}$ in $X_\Gamma$ with the following properties:

%\begin{itemize}
%\item $\I_H(D_i)\searrow b_\Gamma$,
%\item $N\subset \Omega_i$ for any $i>0$,
%\item $Vol(\Omega_i\triangle\Omega_{i+1})\to 0$.
%\end{itemize}

%By using $\Omega_i$, we will construct a special domain $\wh{\Omega}$ in $M$ such that $\Sigma=\partial \wh{\Omega}-S^-$ is desired embedded $H$-disk in $M$ with $\partial \Sigma=\Gamma$.\\

%NOTE THAT WE NEED $H<H_0$ TO CONSTRUCT THE BARRIER N!!! FOR $H=H_0$ ???

\vspace{.2cm}

\noindent {\bf Step 1 - Construction of Lower Barrier $N$ near $S^-$:} In this step, we will show the existence of a neighborhood $N$ of $S^-$ which is away from a sequence of disk $\{D_n\}$ in $Y_\Gamma$ where $\I_H(D_n)\to b_\Gamma$. In other words, we can choose minimizing sequence $\{D_n\}$ with $N\cap D_n=\emptyset$ for any $n>0$.

$S^-$ is a surface in $M$ with boundary $\Gamma$. Consider $S^-$ with the induced metric (and hence the induced topology). Cover $S^-$ with open disks $B_{\e_j}(p_j)$ in $S^-$ where $p_j\in S^-$ and $\e_j>0$ are given by Theorem \ref{H-disk2} for $M$ and $H_0>0$. Since $S^-$ is compact, we can find a finite cover, i.e. $S^-\subset \bigcup_{i=1}^K B_{\e_j}(p_j)$. Let $\gamma_j=\partial B_{\e_j}(p_j)\subset S^-$. Let $\wh{B}_{\e_j}(p_j)$ be the $\e_j$ neighborhood of $p_j$ in $M$. By Theorem \ref{H-disk2}, there is an embedded minimizing $H$-disk $E_j$ in $\wh{B}_{\e_j}(p_j)$ with $\partial E_j = \gamma_j$. By choosing $\{e_j\}$ smaller if necessary, we can assume $E_j$ is also minimizing $H$-disk in $M$ with $\partial E_j=\gamma_j$. Let $N_j$ be the open component of $\wh{B}_{\e_j}(p_j)-E_j$ with $\partial \overline{N_j}\cap S^-=B_{\e_j}(p_j)$. Let $N=\bigcup_{j=1}^K N_j$. We claim that $N$ is the desired barrier near $S^-$.

Now, we show that a given sequence $\{D_i\}\subset Y_\Gamma$ with $\I_H(D_i)\searrow b_\Gamma$ can be modified to be away from $N$. Fix $1\leq j_o\leq K$. Take a sequence  $\{D_i\}\subset Y_\Gamma$ with $\I_H(D_i)\searrow b_\Gamma$. Assume that $D_i\cap N_{j_o}\neq \emptyset$ for some $i>0$. We claim that we can reduce $\I_H(D_i)$ with a simple surgery in $\wh{B}_{\e_{j_o}}(p_{j_o})$.

Consider the minimizing $H$-disk $E_{j_0}\subset \partial \overline{N_{j_0}}$. Take a component $T$ in $N_{j_0}\cap D_i$. Then since $D_i$ is an embedded disk, one of the components in $\partial \overline T\subset E_{j_0}$, say $\alpha$ (the outermost curve in $D_i$), bounds a disk $F$ in $D_i$ such that $T\subset F$. Let $F'$ be the disk in $E_{j_0}$ with boundary $\alpha$. Let $D_i'=(D_i-F)\cup F'$. Since $E_{j_0}$ is a minimizing $H$-disk in $M$, then $\I_H(D_i')\leq I_H(D_i)$. If we repeat the process for all components in $N_{j_0}\cap D_i$, we get an embedded disk $D_i'$ with $\I_H(D_i')\leq I_H(D_i)$. Furthermore, $D_i'\cap N_{j_0}=\emptyset$.

This shows that pushing $D_i$ out of $N_{j_o}$ via surgery does not increase the value of $I_H$. Hence, by doing this surgery for every $j$ for $1\leq j \leq K$, we can obtain a new sequence of embedded disks $\{\wh{D}_i\}$ where $\I_H(\wh{D}_i)\searrow b_\Gamma$ and $\wh{D}_i\cap N=\emptyset$ for any $i>0$. Step 1 follows.\\

%By abuse of notation, we will denote the sequence with $\{D_i\}$ and $\Omega_i$ is defined as before. In particular, we will assume $N\subset \Omega_i$ for any $i>0$ from now on.

\noindent {\bf Step 2 - Construction of Upper Barrier $\wh{\Omega}$:} In this step, we will construct an upper barrier $\wh{\Omega}$ where we can construct a minimizing sequence $\{E_i\}$ in $Y_\Gamma$ with $\I_H(E_i)\searrow b_\Gamma$ and $E_i\subset \wh{\Omega}$ for any $i>0$.

For a given $D_i\in Y_\Gamma$, let $\Omega_i$ be the region in $M$ with $\partial \Omega_i = D_i\cup S^-$ as before. We claim that we can define a sequence of embedded disks $\{D_i\}$ in $Y_\Gamma$ with the following properties.

\begin{itemize}

\item $\I_H(D_i)\searrow b_\Gamma$,

\item $N\subset \Omega_i$ for any $i>0$,

\item $Vol(\Omega_i\triangle\Omega_{i+1})\to 0$.

\end{itemize}

Let $\{D_i\}$ be a sequence in $Y_\Gamma$ with $\I_H(D_i)\searrow b_\Gamma$. By Step 1, we can modify the sequence $\{D_i\}$ by doing surgery, and get a new sequence $\{\wh{D}_i\}$ such that $\I_H(\wh{D}_i)\searrow b_\Gamma$ and $\wh{D}_i\cap N=\emptyset$ for any $i>0$. This ensures the second condition. Finally, for the last condition, we can pass to a subsequence if necessary to guarantee that $Vol(\Omega_i\triangle\Omega_{i+1})\to 0$ where $\triangle$ represents the symmetric difference. Hence, the new sequence, say $\{D_i\}$, has the desired properties.

Now, take such a sequence $\{D_i\}$, and let $\Omega_i$ be as defined before. Let $\wh{\Omega}_n$ be the interior of the component of $\bigcap_{i=n}^\infty \Omega_i$ containing $int(S^-)$. Then by construction, $N\subset \wh{\Omega}_1\subset \wh{\Omega}_2\subset ... \subset \wh{\Omega}_n \subset ...$

Define $\wh{\Omega}=\bigcup_{n=1}^\infty \wh{\Omega}_n$. By construction, $\wh{\Omega}$ is a {\em connected open set} in $M$ with $N\subset \wh{\Omega}$. Notice that again by construction, $Vol(\wh{\Omega}_n\triangle \Omega_n)\to 0$, and $Vol(\wh{\Omega}\triangle \Omega_n)\to 0$ as $n\to \infty$. This shows that we can take a minimizing sequence $\{D_i\}$ in $Y_\Gamma$ where $D_i\subset \wh{\Omega}$ for any $i>0$.

\vspace{.2cm}

%Define $\Sigma=\partial \overline{\wh{\Omega}}-int(S^-)$. We claim that $\Sigma$ is the desired embedded $H$-disk in $M$ with $\partial \Sigma=\Gamma$.\\

%Take a sequence of embedded disks $\{D_i\}\subset X_\Gamma$ with $\I_H(D_i)\searrow b_\Gamma$. Let $\Omega_i$ be the component of $M-D_i$ with $\partial \overline{\Omega}_i= S^-\cup D_i$ for any $i>0$. For the sequence $\{D_i\}$, we will also require that $\Omega_i\supset N$ where $N$ is an open neighborhood of $int(S^-)$ in $M$. In other words, we will show that such a sequence can be chosen such that the disk $D_i$ stays away from $S^-$ by constructing a barrier $N$.

%Define $\wh{\Omega}_n=\bigcap_{i=n}^\infty \overline{\Omega}_i$. Clearly, $\wh{\Omega}_1\subset \wh{\Omega_2}\subset ...\subset \wh{\Omega}_n \subset..$ Let $\wh{\Omega}=\bigcup_{n=1}^\infty \wh{\Omega}_n$.

\noindent {\bf Step 3 - Construction of an ordered sequence of embedded disks $\{E_i\}$ in $Y_\Gamma$:} Now, we will construct a sequence of embedded disks $\{E_i\}$ in $Y_\Gamma$ with the following properties:

\begin{itemize}
\item $\I_H(E_i)\searrow b_\Gamma$.
\item $E_i\cap E_j=\Gamma$ for $i\neq j$.
\item $N\subset\Delta_i\subset\Delta_j\subset \wh{\Omega}$ for $i<j$.
\end{itemize}

Here, $\Delta_i$ is the open component of $M-E_i$ containing $int(S^-)$, i.e. $\partial \overline{\Delta_i}=E_i\cup S^-$. We will call $E_i<E_j$ if $\Delta_i\subset \Delta_j$ and $E_i\cap E_j=\Gamma$.

We will define the sequence inductively. Let $\wt{\Omega}=\overline{\wh{\Omega}}$. Existence of an embedded disk $E_1$ in $\wt{\Omega}-N$ follows from Step 1 and 2. Now, for $1\leq i\leq n$, let $E_i\in Y_\Gamma$ be an embedded disk in $\wt{\Omega}-N$ such that $E_1<E_2<E_3<...<E_n$ and $\I_H(E_{i})<\I_H(E_{i-1})$.

Now, we claim that there is $E_{n+1}\in Y_\Gamma$ such that $E_{n+1}\subset\wt{\Omega}-\Delta_n$ and $\I_H(E_{n+1})<\I_H(E_n)$. For sufficiently large $m>0$, let $D_m\in Y_\Gamma$ be an embedded disk in the construction in Step 2, with $b_\Gamma<\I_H(D_m)<\I_H(E_n)$. If $D_m$ is in $\wt{\Omega}-\Delta_n$, then we are done. If not, notice that by the construction of the open connected region $\wh{\Omega}$, $Vol(\wh{\Omega}\triangle \wh{\Omega}_m)\searrow 0$ as $m\to \infty$. Again, by choosing $m$ sufficiently large, and by pushing $D_m$ slightly into $\wh{\Omega}$ if necessary, we can make sure that the new disk $D_m'$ satisfies the desired properties, i.e. $\I_H(D_m')<\I_H(E_n)$ and $D_m'\subset\wt{\Omega}-\Delta_n$. Hence, define $E_{n+1}=D_m'$. By induction, we get a sequence of embedded disks $E_1<E_2<E_3<...<E_n<..$. Furthermore, we can also assume that $I_H(E_n)<b_\Gamma+\frac{1}{n}$ by choosing $m$ sufficiently large for $D_m$ above. This implies $I_H(E_n)\searrow b_\Gamma$.\\

\noindent {\bf Step 4: The Limit is a surface:} Consider $\wh{\Delta}=\bigcup_{i=1}^\infty \Delta_n$ defined in the previous section. Let $\Sigma=\partial \overline{\wh{\Delta}}- int(S^-)$. In this step, we show that $\Sigma$ is a smooth embedded $H$-surface in $M$ with $\partial \Sigma=\Gamma$. In the next step, we will finish the proof by showing that $\Sigma$ is indeed a disk.

Let $p$ be a point in $\Sigma-\Gamma$. Then, by the construction of $\{E_i\}$, there exists a sequence $\{p_i\}$ such that $p_i\in E_i$ and $d(p_i,p)\searrow 0$. Let $\e_p>0$ be as in Corollary \ref{embed2}. Let $0<\e<\e_p$ sufficiently small so that any minimizing $H$-disk $D$ in $M$ with $\partial D \subset \partial B_\e(p)$ belongs to $B_\e(p)$. Let $E_i^p$  be the component of $E_i\cap B_\e(p)$ containing $p_i$. Then, $E_i^p$ is a compact planar surface in $E_i$. Let $\gamma_i^p\subset \partial B_\e(p)$ be the outermost curve of $\partial E_i^p$ in $E_i$ such that $\gamma_i^p$ bounds a disk $V_i$ in $E_i$ where $E_i^p \subset V_i$.

Now, let $U_i$ be the minimizing $H$-disk in $M$ with $\partial U_i=\gamma_i^p$. By the assumption on $\e$, and Corollary \ref{embed2}, $U_i$ is a properly embedded disk in $B_\e(p)$.

%Furthermore, by taking $\e<\e_p$ smaller if necessary, we can assume that any minimizing $H$-disk $D$ in $M$ with $\partial D\subset \partial B_\e(p)$ belongs to $B_\e(p)$.

Let $E_i'=(E_i-V_i)\cup U_i$. Since $U_i$ is a minimizing $H$-disk in $M$, then $\I_H(E_i')<I_H(E_i)$. Furthermore, since $\gamma_i^p<\gamma_j^p$ for $i<j$ in $\partial B_\e(p)$ by the order construction, simple comparison shows that $U_i<U_j$ in $B_\e(p)$ likewise. Hence, with this modification, the order has been kept in the sequence of disks $\{E_i'\}$, i.e. $\Delta_i'\subset \Delta_j'$ for $i<j$. Furthermore, as $\I_H(E_i')\leq I_H(E_i)$, $I_H(E_i')\searrow b_\Gamma$. Hence, this proves that $\wh{\Delta}=\bigcup_{i=1}^\infty \Delta_i'$. Then, we can consider $\Sigma\cap B_\e(p)$ as the limit of the minimizing $H$-disks $U_i=E_i'\cap B_\e(p)$. Hence, if we show that the limit of minimizing $H$-disks $\{U_i\}$ is smooth embedded $H$-disk, this implies $\Sigma$ is a smoothly embedded surface in $M$.

We will adapt the proof of \cite[Lemma 3.3]{HS} to $H$-disks where Hass and Scott proved the similar statement for a sequence of least area disks.

%Recall that $\partial U_i=\gamma_i^p$ is a simple closed curve in $\partial \B_\e(p)$, and $\partial \B_\e(p)-\gamma_i^p=Y_i^-\cup Y_i^+$ where $Y_i^\pm$ is an open disk in $\partial \B_\e(p)$. By construction $U_i<U_j$, $Y_i^-\subset Y_j^-$ for any $i<j$. Hence, $S_1^-\subset S_2^-\subset ... \subset S_i^-\subset...$. Then, $S^-=\bigcup_{i=1}^\infty S_i^-$ is an open disk in $\partial \B_\e(p)$. Let $\wh{\gamma}=\partial \overline{S^-}$ in $\partial \B_\e(p)$. If $\wh{\gamma}$ is a point, then we say the sequence $\{D\}$ escapes out of $\B_\e{p}$, and the limit is empty. Otherwise, there is a limit point $p$ in $int(\B_\e(p))$ of the sequence $\{D_i\}$, i.e. $\exists p_{i_j}\in D_{i_j}$ with $p_{i_j}\to p$.

By construction $p$ is also a limit point of $\{U_i\}$ as $\wh{\Delta}=\bigcup_{i=1}^\infty \Delta_i'$ and $p\in \Sigma\subset \partial \overline{\wh{\Delta}}$ i.e. $\exists \{q_i \ | \ q_i\in U_i\}$ such that $q_i\to p$. Since $\{U_i\}$ is a sequence of minimizing $H$-disks, for any $U_i$, we have a conformal $H$-harmonic parametrization $u_i:D^2\to M$ with $u_i(D^2)=U_i$. By composing with a conformal automorphism of disk, we can assume $u_i(0)=q_i\to p$. As $u_i$ conformal, the energy and the area are equal. Then as $\{U_i\}$ are minimizing $H$-disks, there is a uniform bound $K$ such that $E(u_i)<K$ where $E(.)$ is the energy. Then, as in \cite[Lemma 3.3]{HS}, by using Courant-Lebesque lemma, the equicontinuity of the family $\{u_i\}$ on the compact subsets of $int(D^2)$ follows. Hence, by Arzela-Ascoli theorem, there is a continuous function $u:int(D^2)\to M$ where $u_i\to u$ pointwise on $int(D^2)$. Furthermore, by using the fact that $\{u_i\}$ are conformal $H$-harmonic maps, the arguments in \cite{HS} shows that $\{u_i\}$ are bounded in the $C^{2,\alpha}$ norm, and $\{u_i\}$ has a subsequence converging $u$ smoothly on $int(D)$.

Furthermore, as $\{U_i\}$ are minimizing embedded $H$-disks in $M$, the limit $u(D^2)=\wh{U}$ is also minimizing embedded $H$-disk by simple comparison again as in \cite[Lemma 3.3]{HS}. Let $\partial \B_\e(p)-\gamma_i^p=Y_i^-\cup Y_i^+$ where $Y_i^\pm$ is an open disk in $\partial B_\e(p)$ with $\partial \overline{\Delta_i'\cap B_\e(p)}=Y_i^-$. By the definition of the ordered sequence $\{E_i\}$, $U_i<U_j$ in $B_\e(p)$, and $Y_i^-\subset Y_j^-$ for any $i<j$. Hence $Y_1^-\subset Y_2^-\subset...\subset Y_i^-\subset..$, and let $\wh{Y}=\bigcup_{i=1}^\infty Y_i^-$. Then, $\wh{Y}$ is an open disk in $\partial B_\e(p)$ with $\partial \wh{U}=\partial \overline{\wh{Y}}$.

Let $\wh{\gamma}=\partial \overline{\wh{Y}}$. Since $Y$ is an open disk in $\partial B_\e(p)$, $\wh{\gamma}$ is a closed curve, which may not be simple (even it might be a point). Since $p\not\in \wh{\Delta}$, $Y_i^-=\partial \wh{\Delta}$ cannot be the whole $\partial B_\e(p)$ except a point. Hence, $\wh{\gamma}$ is not a single point. If $\wh{\gamma}$ is not a {\em simple} closed curve, we can pass to a smaller $\e>0$, and make sure that $\wh{\gamma}$ is a simple closed curve as $\wh{U}$ is the limit of embedded minimizing $H$-disks. This proves that $\Sigma\cap B_\e(p)=\wh{U}_p$ is a smooth embedded disk. Hence, for any $p\in int(\Sigma)$ we can find such neighborhood $\wh{U}_p$ in $\Sigma$. This shows that $\Sigma$ is a smoothly embedded surface in the interior. $\Sigma$ is a continuous surface with boundary $\Gamma$ by construction as for any embedded disk in the ordered sequence of embedded disks $\{E_i\}$ in $M$ with $\partial E_i=\Gamma$. Hence, $\Sigma$ is a properly embedded smooth $H$-surface in $M$ with $\partial \Sigma=\Gamma$.\\

\noindent {\bf Step 5 - Embedded $H$-disk in $M$:} Now, we will show that $\Sigma$ is indeed a disk, and finish the proof of the theorem

In the previous step, we have shown that $\Sigma$ is a properly embedded smooth $H$-surface in $M$ with $\partial \Sigma=\Gamma$. If $\Sigma$ is not a disk, then it must be a genus $g$ surface with $1$ boundary component where $g>0$. Furthermore, $\Sigma\subset \partial \overline{\wh{\Delta}}$. Hence, there must be a smooth simple closed curve $\alpha$ in $\Sigma$ which does not bound any disk in $\overline{\wh{\Delta}}$ as $g>0$.

Now, let $\e_o>0$ sufficiently small so that the $\e_o$ neighborhood $N_{\e_o}(\alpha)$ of $\alpha$ in $M$ is a solid torus. By construction, for sufficiently large $n_o>0$, there are embedded disks $E_{n_o}\subset \wh{\Delta}$ with $E_{n_o}\cap N_{\e_o}(\alpha)$ contains a simple closed curve $\alpha_o$ close to $\alpha$. Furthermore, we can make $\alpha_o\subset \wh{\Delta}$ as close as we want by choosing $\e_o>0$ small enough. Hence, there is an annulus $\mathcal{A}$ in $N_{\e_o}(\alpha)\cap\overline{\wh{\Delta}}$ with $\partial \mathcal{A}=\alpha\cup\alpha_o$. Furthermore, since $E_{n_o}$ is an embedded disk in $M$, $\alpha_o$ bounds a disk $D_o\subset E_{n_o}$ with $\partial D_o=\alpha_o$. Hence, $\wh{D}=\mathcal{A}_o\cup D_o$ is a disk in $\overline{\wh{\Delta}}$ with $\partial \wh{D}=\alpha$. This is a contradiction. Hence, $\Sigma$ must be a disk. The proof follows.
\end{pf}

\begin{rmk} [Existence of solutions to $H$-Plateau problem in $H_0$-convex domains] Notice that the results in Section 3 shows the embeddedness of the immersed solutions of $H$-Plateau problem, whose existence obtained in the classical papers of the subject. However, in the result above, we not only show embeddedness, but also show the existence of such solutions for the first time in this setting.

In particular, before this theorem, there was no result on the existence of $H$-disks for such general domains, even for the immersed case. The most general results on the existence of immersed $H$-disks are the round balls in $\BR^3$ (Lemma \ref{H-disk1} \cite{Gu}), and small balls in Riemannian $3$-manifolds where the size is determined by curvature of the ambient space (Lemma \ref{H-disk2} \cite{Gu, HK}). In this respect, our result only assumes a global convexity condition on the boundary $\partial M$, and gives a positive answer to the $H$-Plateau problem for the nullhomotopic $H_0$-extreme curves, i.e. the curves in $\partial M$. Furthermore, the solutions obtained are embedded.
\end{rmk}

\subsection{Rellich Conjecture} \label{Rellich-sec} \

\vspace{.2cm}

Rellich conjectured that any simple closed curve $\Gamma$ bounds at least two solutions of $H$-Plateau problem for $H$ sufficiently small \cite{BC}. Our techniques naturally provides a positive answer to this conjecture in more general settings. Furthermore, the pair of $H$-disks are both embedded.

\begin{cor} \label{rellich} Let $M$ be an $H_0$-convex manifold with $H_2(M)=\{0\}$. Let $\Gamma$ be a Jordan curve in $\partial M$ such that $\Gamma$ is nullhomotopic in $M$ and separating in $\partial M$. Then, for any $H\in (0,H_0)$, there are two embedded $H$-disks $\Sigma^\pm_H$ in $M$ with $\partial \Sigma^\pm_H=\Gamma$.
\end{cor}

\begin{pf} In the proof of Theorem \ref{main}, in the construction of $\wh{\Omega}$, if we use $S^+$ instead of $S^-$ where $S-\Gamma=S^+\cup S^-$, then, we would get a different embedded $H$-disk $\Sigma^+$ with the same boundary $\Gamma$ for any $H\in (0,H_0)$. Hence, if we call the $H$-disk obtained in Theorem \ref{main} as $\Sigma^-$, then both $\Sigma^+$ and $\Sigma^-$ would be embedded $H$-disks in $M$ with $\partial \Sigma^\pm=\Gamma$ where the convex sides facing each other. For $H>0$, $\Sigma^+$ and $\Sigma^-$ cannot be the same disk as the $H$-convex sides are pointing opposite directions. Notice that for $H=0$, they could be the same disk, which in that case, $\Gamma$ bounds a unique area minimizing disk in $M$ \cite{Co}.
\end{pf}

\subsection{Dehn's Lemma and the nonexistence of minimizing $H$-disks} \label{Dehn} \

\vspace{.2cm}

In this part, we will show the that Theorem \ref{main} is indeed {\em sharp} in the sense that the homology conditions on $M$ and $\Gamma$ in the statement cannot be removed to prove the Dehn's Lemma for $H$-disks in full generality (See Remark \ref{dehn3}). We first recall the  Dehn's Lemma, which is one of the fundamental results in $3$-manifold topology.

\vspace{.2cm}

\noindent {\em Dehn's Lemma:} \cite{Pa} Let $M$ be a $3$-manifold with boundary. Let $\gamma$ be a simple closed curve in $\partial M$. If $\gamma$ is nullhomotopic in $M$, then $\gamma$ bounds an embedded disk in $M$.

\vspace{.2cm}

Notice that the conditions $\Gamma$ is separating in $\partial M$, and $H_2(M)=\{0\}$ in Theorem \ref{main} do not appear the original statement of Dehn's Lemma. Now, we will show that these two conditions are necessary to prove the {\em existence} of minimizing $H$-disk $\Sigma_H$ in an $H_0$-convex manifold $M$ with $\partial \Sigma_H=\Gamma$. In particular, we will construct a counterexample as follows.

\vspace{.2cm}

\noindent {\bf Claim:} There exists a compact $2$-convex $3$-manifold $\M$ and a Jordan curve $\gamma$ in $\partial \M$ where $\gamma$ is nullhomotopic in $\M$ such that there exists no minimizing $H$-disk $\Sigma_H$ in $\M$ with $\partial \Sigma_H=\gamma$ for $H\in(1,2)$.

\vspace{.2cm}

\noindent {\em Proof.} Let $\B_r$ be the closed $r$-disk in $\BR^2$ with center origin. Let $M=\B_2\times \BR$ be an infinite solid cylinder of radius $2$ in $\BR^3$. Hence, $\partial M$ is the infinite cylinder in $\BR^3$ with radius $2$. Let $\Gamma$ be the round circle $\partial \B_2\times \{0\}$.

Now, we will modify the metric on $M$ to define a $2$-convex manifold $\wh{M}$. Consider $M'=\B_\frac{1}{4}\times \BR$. Clearly, $M'$ is a $2$-convex manifold. Furthermore, both $M$ and $M'$ are both rotationally invariant, and vertical translation invariant. Let $\e>0$ be sufficiently small. Consider the $2\e$ neighborhood $N$ of $\partial M$ in $M$, and the $2\e$ neighborhood $N'$ of $\partial M'$ in $M$, i.e. $N=M-\B_{2-2\e}\times \BR$ and $N'=M'-\B_{\frac{1}{4}-2\e}\times \BR$. Use cylindrical coordinates in both manifolds. Let $\varphi:N\to N'$ be a diffeomorphism such that $\varphi(r, \theta, t)=(r-\frac{7}{4},\theta,t)$. Let $g'$ be the metric on $N$ induced by $\varphi$.

Now, we will keep the original metric $g$ on $\B_{2-2\e}\times \BR$, and the new metric $g'$ on $N$, but we will modify the metric on the region between them, on $(\B_{2-\e}-\B_{2-2\e})\times \BR$, by using a partition of unity so that the modified metric on $M$ would be smooth. Define a partition of unity, $\psi_1$ and $\psi_2$, smooth positive functions on $[0,2]$ with $\psi_1(r)+\psi_2(r)=1$ as follows.
%$\psi_1$ and $\psi_2$ are smooth functions on $[0,2]$ with $0\leq \psi_1(r),\psi_2(r)\leq 1$ and $\psi_1(r)+\psi_2(r)=1$ for any $r\in [0,2]$.
%Let $\psi_1(r)=1$ on $[0,2-2\e]$ and $\psi_1(r)=0$ on $[2-\e,2]$. Similarly, let $\psi_2(r)=0$ on $[0,2-2\e]$ and $\psi_2(r)=1$ on $[2-\e,2]$. In particular,

$$\ \ \ \psi_1(r)=\left\{ \begin{array}{ll}
1 \ \  & r\in[0,2-2\e]\\
0 \ \  & r\in[2-\e,2]\\
\end{array}\right .\ \ \ \
\psi_2(r)=\left\{ \begin{array}{ll}
0 \ \  & r\in[0,2-2\e]\\
1 \ \  & r\in[2-\e,2]\\
\end{array}\right .$$

Now define the smooth metric $\wh{g}$ on $M$ such that $\wh{g}=\psi_1.g+\psi_2.g'$ where $g$ is the original product metric on $M$. Let $\wh{M}$ be the manifold $M$ with this new metric $\wh{g}$. By construction, $\wh{M}$ is $2$-convex. Furthermore, $\wh{M}$ is both rotationally invariant, and vertical translation invariant by construction. In particular, $\phi_{\theta_o}(r,\theta, t)=(r, \theta+\theta_o, t)$ and $F_c(r,\theta,t)=(r,\theta, t+C)$ are both isometries of $\wh{M}$ for any $\theta_o,C\in \BR$.

We claim that there is no minimizing $H$-disk $\Sigma_H$ in $\wh{M}$ with $\partial \Sigma_H=\Gamma$ for $H\in(1,2)$ even though $\wh{M}$ is a $2$-convex manifold.
Now, assume on the contrary that there exists a minimizing $H$-disk $\Sigma_H$ in $\wh{M}$ with $\partial \Sigma_H=\Gamma$ for $H\in(1,2)$. Let $X_\Gamma$ be the space of the immersed disks in $\wh{M}$ with boundary $\Gamma$.

Hence, by definition, $\I_H(\Sigma_H)$ achieves the infimum in $X_\Gamma$. Notice that to define $\I_H$, we first need to fix a surface $T$ in $\wh{M}$ with $\partial T=\Gamma$. Hence, for $D\in X_\Gamma$, define $\Delta_D$ be the region in $\wh{M}$ with $\partial \Delta_D=T\cup D$. As $H_2(M)=\{0\}$, $T\cup D$ always separates a region $\Delta_D$ in $\wh{M}$. Hence, $\I_H(D)=|D|+2H\|\Delta_D\|$ for any $D\in X_\Gamma$. Notice that the minimizer of $\I_H$ is independent of the choice of the surface $T$, as for some other choice, the functional $\I_H$ would change by a constant.

As $\wh{M}$ is mean convex and homogeneously regular, there exists an area minimizing disk $\V$ in $\wh{M}$ by \cite{MY2}. In particular, it is easy to see that $\V=\B_2\times \{0\}$ by the construction of $\wh{M}$. Without loss of generality, we can assume $T\cap \V=\Gamma$. Let $\wh{\Delta}$ be the region in $\wh{M}$ with $\partial\wh{\Delta}=\V\cup T$. Let $C_0=\|\wh{\Delta}\|$, the volume of $\wh{\Delta}$. By a simple area comparison, we can assume that minimizers of $\I_H$ for $H>0$ stays in one side of $\V$. Hence, for any $D\in X_\Gamma$, we restrict ourselves to $D\subset \wh{\Delta}$, and hence $\Delta_D\subset \wh{\Delta}$. For such $D\in X_\Gamma$, let $\wh{\Delta}_D=\wh{\Delta}-\Delta_D$, the region between $D$ and $\V$. In particular, $\partial \overline{\wh{\Delta}_D}=D\cup\V$ and $\|\wh{\Delta}_D\|=C_0-\|\Delta_D\|$. Now, consider $$\I_H(D)=|D|+2H\|\Delta_D\|=|D|+2H(C_0-\|\wh{\Delta}_D\|)=2HC_0+(|D|-2H\|\wh{\Delta}_D\|)$$

Let $\wh{\I}_H(D)=|D|-2H\|\wh{\Delta}_D\|$. Then, $\I_H(D)=2HC_0+\wh{\I}_H(D)$. Hence,  $\Sigma_H$ minimizes $\I_H$ if and only if it minimizes $\wh{\I}_H$. Now, we show that there is a sequence $E_n\in X_\Gamma$ where $\wh{\I}_H(E_n)\searrow -\infty$, and get a contradiction. This will prove that there is no minimizing $H$-disk in $\wh{M}$ with boundary $\Gamma$.

Let $A=\V-int(\B_1\times\{0\})$ be the annulus in $\V$. Let $\C_n=\B_1\times [-n,0]$ be the solid cylinder in $\wh{M}$ of radius $1$ and height $n$. Let $E_n'=\partial \C_n - int(\B_1\times\{0\})$. Define $E_n=A\cup E_n'$. Consider $\wh{\I}_H(E_n)= |E_n|-2H\|\Delta_n\|=|E_n|-2H\|\C_n\|$ where $\Delta_n=\wh{\Delta}_{E_n}=\C_n$ is the region between $\V$ and $E_n$. By straightforward calculation, $|E_n|=c_o+2\pi n$ where $c_0=|A|+|\B_1\times\{-n\}|$, and $2\pi n$ is the area of the cylinder radius $1$ and of height $n$, i.e $\partial \C_n$ with no caps. Also, $\|\Delta_n\|=\|\C_n\|=\pi n$, the volume of $\C_n$. Hence, $$\wh{\I}_H(E_n)= |E_n|-2H\|\Delta_n\|=(c_o+2\pi n)- 2H\pi n=c_0+2\pi n(1-H)$$ This implies that $\wh{\I}_H(E_n)\searrow-\infty$ for $H\in(1,2)$.

Without loss of generality, by smoothing out the nonsmooth parts in $E_n$, we can assume $E_n\in X_\Gamma$. Hence, $\inf_{X_\Gamma}\wh{\I}_H=-\infty$. Now, if we recall the relation between $\I_H$ and $\F_H$ in section \ref{H-Plateau}, we see that $\F_H(E_n)\searrow -\infty$ for $H\in(1,2)$, too. The advantage of $\F_H$ over $\I_H$ is that the functional $\F_H$ does not require a surface $T$ to define the region $\Delta_D$ used in $\I_H$. On the other hand, if $\Sigma_H$ is a minimizing $H$-disk in $\wh{M}$, it must achieve the $\inf_{X_\Gamma}\F_H$ for a conformal parametrization of $\Sigma_H$ by Lemma \ref{variation}. This shows that there is no minimizing $H$-disk $\Sigma_H$ in $\wh{M}$ with $\Sigma_H=\Gamma$. %Recall that for any choice of initial surface $T$ and $T'$, $\I_H-\I_H'$ is constant, which is the volume of the region between $T$ and $T'$.

Now, $\wh{M}$ is not a compact $3$-manifold. By taking a quotient of $\wh{M}$ by the vertical $\BZ$-action, we get a compact manifold $\M$ with the desired properties, and finish the claim. Let $F_n:\wh{M}\to\wh{M}$ be the isometry of $\wh{M}$ corresponding vertical translation by $n\in \BZ$, i.e. $F_n(r,\theta,t)=(r,\theta, t+n)$. Let $G$ be the group of integer vertical translations, i.e. $G=<F_n \ | \ n\in\BZ>\simeq\BZ$. Consider the compact manifold $\M=\wh{M}/G$. Hence, $\M$ is a compact $3$-manifold, which is topologically a solid torus.  Furthermore, $\M$ is a $2$-convex manifold.

By construction, $\wh{M}$ is the universal cover of $\M$, and let $\pi:\wh{M}\to \M$ be the covering map. Let $\gamma=\pi(\Gamma)$ be the simple closed curve in $\partial \M$. Like before, we claim that there is no minimizing $H$-disk in $\M$ with boundary $\gamma$ for $H\in(1,2)$. In the construction above, we can modify $E_n$ so that $\pi(E_n)$ is an embedded disk in $\M$. In particular, for sufficiently small $\delta>0$, define the slanted cylinder $\C_n^\delta$ with a generating line $l_n^\delta$ connecting points $(1-\delta, 0, -n+\delta)$ and $(1+\delta, 0, 0)$ in the $\theta=0$ slice in $\wh{M}$. Hence, $\C_n^\delta$ is the slanted cylinder obtained by rotating $l_n^\delta$ around the $z$-axis. If we replace $\C_n^\delta$ with $\C_n$ in the definition of $E_n$, we get $F_i(E_n)\cap F_j(E_n)=\emptyset$ for any $i\neq j\in \BZ$. This implies $\pi(E_n)=D_n$ is embedded disk in $\M$ with $\partial D_n=\gamma$. Furthermore, $\F_H(D_n)=\F_H(E_n)\searrow -\infty$ for $H\in(1,2)$ as before. This shows that there exists no minimizing $H$-disk $\Sigma_H$ in $\M$ with $\partial \Sigma_H=\gamma$ for $H\in(1,2)$. The claim follows. 

\hfill $\Box$

\begin{rmk} \label{dehn3} In \cite{MY1,MY2,MY3}, Meeks and Yau proved Dehn's Lemma for area minimizing disks in mean convex manifolds. The existence of area minimizing disks is given by Morrey's famous result \cite{Mo}. Meeks and Yau proved that Morrey's solutions for the Plateau problem are embedded.

In Theorem \ref{main}, we generalize both of these results. Notice that Theorem \ref{main} proves the existence of minimizing $H$-disks in $H_0$-c0nvex manifolds, and embeddedness of these results follows from Theorem \ref{embed1}. 

In order to show the existence, the homology conditions $\Gamma$ being separating in $\partial M$, and $H_2(M)=\{0\}$ are essential. This is because to minimize $\I_H$, these conditions induces a lower barrier, and automatically defines a lower bound for $\I_H$, i.e. $\I_H(D)\geq 0$. The example above shows that if we remove these conditions, the existence of the minimizing $H$-disk cannot be guaranteed. e.g. in the example above, $\M$ is a solid torus, and $\Gamma$ is the meridian curve in $\partial \M$, which is not separating in $\partial \M$. See section \ref{Dehn2} for further discussion.
\end{rmk}

\section{Concluding Remarks}

\subsection{Embedded Planar $H$-surfaces} \

By using similar arguments, the main theorem (Theorem \ref{main}) can be naturally generalized to planar domains for a given collection of pairwise disjoint simple closed curves $\Gamma=\gamma_1\cup \gamma_2\cup ...\gamma_n$ in the boundary of $H_0$-convex domain $M$.  To generalize it to the planar domains, the essential condition would be the infimum of $I_H$ taken over all planar surfaces (connected or not connected) should be achieved for a connected planar surface as in \cite[Theorem 1]{MY3}.

In particular, let $M$ be a $H_0$-convex manifold with $H_2(M)=\{0\}$. Let $\Gamma=\gamma_1\cup \gamma_2\cup ...\gamma_n$ be a collection of pairwise disjoint Jordan curves in $\partial M$ such that $\Gamma$ is separating in $\partial M$, and $\Gamma$ bounds a connected embedded planar surface in $M$. Let $X_\Gamma$ be the space of all immersed planar surfaces with boundary $\Gamma$. Let $Y_\Gamma$ be the space of all immersed planar surfaces with boundary $\Gamma$ where the domain of the immersion is compact, connected planar surface $\mathcal{P}_n$ with $n$ boundary components, i.e. $\mathcal{P}_n$ is a closed disk with $n-1$ open disjoint disks removed.

Hence, the generalization of \cite[Theorem 1]{MY3} to CMC setting should be as follows: If $\inf_{P\in Y_\Gamma}\I_H(P)<\inf_{P\in X_\Gamma-Y_\Gamma}\I_H(P)$, then for any $0\leq H<H_0$, any minimizing planar $H$-surface $\Sigma_H$ in $M$ with $\partial \Sigma_H=\Gamma$ is embedded and $\Sigma_H\in Y_\Gamma$.

The idea is again to construct a lower barrier $N$, and upper barrier $\wh{\Omega}$ as in Theorem \ref{main} by using a sequence of embedded connected planar surfaces $P_i$ with  $\partial P_i=\Gamma$ and $\I_H(P_i)\to inf_{Y_\Gamma}\I_H(P)$. By using $\inf_{P\in Y_\Gamma}\I_H(P)<\inf_{P\in X_\Gamma-Y_\Gamma}\I_H(P)$ condition as in \cite[Theorem 1]{MY3}, the pinching of the sequence can be ruled out. Then with similar arguments, the embeddedness of the limit can be shown  by using the fact that the surfaces are planar.

\vspace{-.2cm}

\subsection{Generalization of the Dehn's Lemma} \label{Dehn2} \

%\vspace{.2cm}

Notice that in section \ref{Dehn}, we showed that Theorem \ref{main} is sharp in some sense for generalization of the Dehn's Lemma. In particular, we constructed an example of an $H_0$-convex manifold $M$, and a Jordan curve in $\partial M$, nullhomotopic in $M$, which does not bound any \underline{minimizing} $H$-disk in $M$ for some values $H<H_0$. However, this does not completely prove that the generalization of Dehn's Lemma for $H$-disks is not possible. In that example, the nonexistence of minimizing $H$-disks was shown, but there still might be {\em non-minimizing} embedded $H$-disks for the given boundary curve.

By Gulliver's result, $H_0$-convexity condition is necessary for the ambient manifold, as for $H>1$, there is no $H$-disk in the unit ball $\B$ in $\BR^3$ bounding the equator circle in $\partial \B$ \cite{Gu}. Hence, the validity of the following version of Dehn's Lemma is a very interesting question.

\vspace{.2cm}

\noindent {\em Dehn's Lemma for $H$-disks:} Let $M$ be a compact $H_0$-convex manifold. Let $\gamma$ be a Jordan curve in $\partial M$. If $\gamma$ is nullhomotopic in $M$, then  for any $H\in[0,H_0)$, there exists an embedded $H$-disk $\D_H$ in $M$ with $\partial \D_H=\gamma$.

%We believe that our main result (Theorem \ref{main}) can be generalized to give a positive answer to Dehn's Lemma for $H$-disks in $H_0$-convex manifolds. In particular, Dehn's Lemma states that if $\Gamma$ is a Jordan curve in the boundary of a $3$-manifold $M$, and if $\Gamma$ is nullhomotopic in $M$, then $\Gamma$ bounds an embedded disk in $M$. Notice that the condition $\Gamma$ is separating in $\partial M$, and the condition $H_2(M)=\{0\}$ in our main result do not appear the original Dehn's Lemma where Meeks and Yau proved for area minimizing surfaces \cite{MY3}. However, in our construction we needed these conditions to get the domain $\Omega_i$ and define $\I_H$ for $D_i$ in the proof above. We believe to employ these techniques the condition $H_2(M)=\{0\}$ is essential, but the condition $\Gamma$ is separating in $\partial M$ might be bypassed by using a homologous curve $\Gamma'$ in $\partial M$ "far away" from $\Gamma$ so that the domain $\Omega_i$ can be defined. Furthermore, in order to get the full generality, $H_2(M)=\{0\}$ condition might be bypassed by using some other techniques, like using $\F_H$ instead of $\I_H$.

\vspace{-.2cm}

\subsection{Other Generalizations} \

%\vspace{.2cm}

The proof of the main result shows that our techniques can be useful to construct minimizing $H$-surfaces in homotopy classes. Hence, by generalizing Hass and Scott's techniques \cite{HS} to this context in a suitable way can be used to construct different special $H$-surfaces in $H_0$-convex manifolds.

For example, for sufficiently small $H>0$ depending on the ambient Haken manifold, it might be possible to construct incompressible $H$-surfaces by constructing a minimizer for the $\I_H$ in a suitable setting. The main problem here would be to define $\I_H$ as there is no domain $\Omega$ to be used. However, as in the construction of the counterexample in Section \ref{Dehn}, by fixing the area minimizing surface $S$ in one side, and solving the variational problem $\wh{\I}_H(\Sigma)= Area(\Sigma)-2HVol(\Omega)$ where $\Omega$ is the region between $\Sigma$ and $S$, it might be possible to construct an $H$-surface in the homology class of $S$. However, one first needs to generalize Lemma 3.3 and Lemma 4.1 in \cite{HS} to this setting.

On the other hand, it might also be possible to show the existence of an $H$-surface in a given homotopy class by using again Hass and Scott's techniques \cite[Theorem 6.13]{HS}. In particular, let $M$ be an $H_0$-convex manifold, and let $\Gamma$ be a collection of simple closed curves in $\partial M$. Let $X_\Gamma=\{u:\Sigma\to M \ | \ u(\partial \Sigma )=\Gamma\}$ be the space of immersions (or homotopy class of a given map $u_0$). Then, it might be possible to show the existence of an immersion which minimizes $\I_H$ in $X_\Gamma$.

\subsection{Nonexamples} \label{nonexample} \

After generalizing Douglas and Morrey's results on original Plateau problem to the $H$-Plateau problem by showing the existence of immersed $H$-disks for suitable domains, the question of the embeddedness of these solutions was mostly ignored because a generalization of Meeks and Yau's embeddedness results was thought to be not possible for $H$-Plateau setting. The main reason for that it was believed that some well-known Jordan curves in $\partial \B$ do not bound any embedded $H$-disk in $\B$ for $H>0$. Here, the misleading point is that for the curves mentioned, the "obvious" $H$-disks bounding these curves are not embedded. However, these disks are not the only $H$-disks they bound, and the minimizing $H$-disks are indeed embedded. Hence, our results might be considered surprising in the field.

\begin{figure}[b]
	\begin{center}
		$\begin{array}{c@{\hspace{.4in}}c}
		
		\relabelbox  {\epsfysize=1.8in \epsfbox{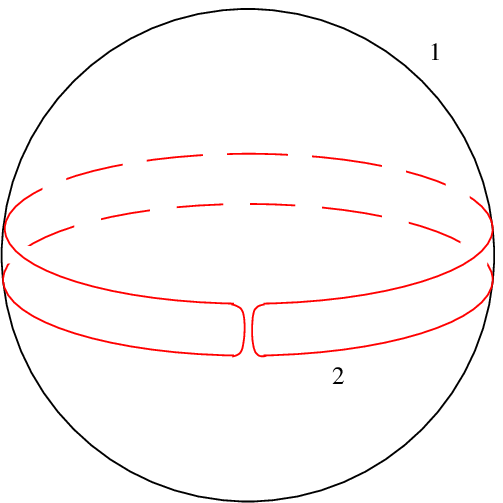}} \relabel{1}{\tiny $\B$} \relabel{2}{\tiny $\Gamma_1$}  \endrelabelbox &
		
		\relabelbox  {\epsfysize=1.8in \epsfbox{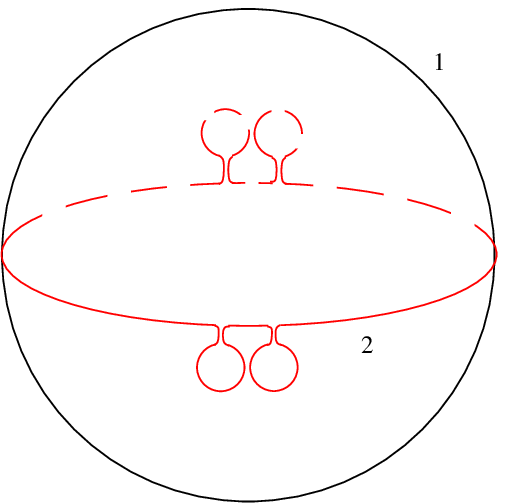}} \relabel{1}{\tiny $\B$} \relabel{2}{\tiny $\Gamma_2$} \endrelabelbox \\
		\end{array}$
		
	\end{center}
	\caption{\label{nonexamples} \footnotesize $\Gamma_1$ and $\Gamma_2$ are Jordan curves in the boundary of the unit ball $\B$. For $H\in(0,1)$, the obvious examples of $H$-disks $\Sigma_i$ in $\B$ with $\partial \Sigma_i=\Gamma_i$ are not embedded. However, the minimizing $H$-disks are always embedded by Theorem \ref{main}.}
\end{figure}

The figures below illustrate some examples of such curves which were believed to have no embedded solution for $H$-Plateau problem. In Figure \ref{nonexamples}, both $\Gamma_1$ and $\Gamma_2$ are Jordan curves in the boundary of the unit ball $\B$.

\vspace{.2cm}

\noindent {\em Example 1:} Let $\gamma^+$ and $\gamma^-$ be two round horizontal circles in $\partial \B$ very close to the equator. $\Gamma_1$ is obtained by connecting  $\gamma^+$ and $\gamma^-$ in $\partial \B$ via a bridge. Fix $H\in(0,1)$. $\gamma^+$ and $\gamma^-$ bounds $H$-disks $D^+$ and $D^-$ which are spherical caps of radius $\frac{1}{H}$ where $\mathbf{H}$ points opposite directions.  We can choose $\gamma^+$ and $\gamma^-$ sufficiently close so that $D^+\cap D^-\neq \emptyset$. Then, by adding a bridge connecting $D^+$ and $D^-$, we can get a nonembedded $H$-disk $\Sigma_1$ in $\B$ with $\Sigma=\Gamma_1$. Notice that we had to choose $D^+$ and $D^-$ for $\gamma^+$ and $\gamma^-$ so that $\mathbf{H}$ points opposite directions to match orientation on $\Sigma_1$ after adding the bridge.

On the other hand, the minimizing $H$-disks bounding $\Gamma_1$ are both embedded. If $\partial \B-\Gamma_1=\Delta^+\cup \Delta^-$, then depending on the choice, the minimizing $H$-disk $D_H$ with $\partial D_H=\Gamma_1$ is either a graph over $\Delta^-$ or a graph over $\Delta^+$. Furthermore, if $E$ is the area minimizing disk in $\B$ with $\partial E=\Gamma$, then these embedded $H$-minimizing disks are in the opposite sides of $E$ by Corollary \ref{rellich}. There might be a nonembedded immersed $H$-disks like $\Sigma_1$ with $\partial \Sigma_1=\Gamma_1$, but it does not imply the nonexistence of embedded $H$-disks $D$ in $\B$ with $\partial D=\Gamma_1$.

\vspace{.2cm}

\noindent {\em Example 2:} The second example is basically more delicate version of the first one. In the first example, because of the asymmetry, the existence of embedded $H$-disk can be seen by considering the larger domain, say $\Delta^+$, where $\partial \B-\Gamma_1=\Delta^+\cup \Delta^-$. Hence, by attaching the bridges to {\em other} $H$-disks $\gamma^+$ and $\gamma^-$ bounds which are close to $\Delta^+$, the constructed $H$-disk $\Sigma_1'$ would be embedded. However, in the second example, we have a symmetric picture, and by choosing $\Omega^+$ or $\Omega^-$ makes no difference where $\partial \B-\Gamma_2=\Omega^+\cup \Omega^-$. Let $\gamma_1^\pm$ and $\gamma_2^\pm$ two pairs of circles in the opposite side of the equator circle $\tau$ where $\gamma_i^+$ are the pair in the positive side, and $\gamma_i^-$ are the pair in the negative side. Then, we obtain $\Gamma_2$ by attaching these circles to $\tau$ as in the figure by keeping the symmetry. Notice that as $\gamma_i^\pm$ is a round circle, the $H$-disks $D_i^\pm$ bounding $\gamma_i^\pm$ are either convex or concave spherical caps depending on the choice. Now, for given $H\in(0,1)$, we can choose $\gamma_1^+$ sufficiently close to $\gamma_2^+$ (and hence $\gamma_1^-$ sufficiently close to $\gamma_2^-$), so that either $D_1^+\cap D_2^+\neq \emptyset$ or $D_1^-\cap D_2^-\neq \emptyset$. Hence, in the second example, in both choices $\Omega^+$ and $\Omega^-$, it seems like there is no embedded $H$-disk in $\B$ with boundary $\Gamma_2$.

\end{document}